\documentclass[onefignum,onetabnum]{siamart220329}

\usepackage[mode=multiuser,status=draft]{fixme} 
\fxsetup{innerlayout=inline}
\fxsetup{targetlayout=colorcb}
\fxusetheme{color}
\FXRegisterAuthor{cc}{envcc}{CC}
\FXRegisterAuthor{gk}{envgk}{GK}

\usepackage{braket,amsfonts}

\usepackage{array}

\usepackage{mathtools}

\usepackage{arydshln}

\usepackage{multirow}

\usepackage{listings}

\usepackage{tikz,tikzscale}
\usepackage{graphicx}

\usetikzlibrary{calc}
\usetikzlibrary{shadings}
\usetikzlibrary{shapes.geometric}
\usetikzlibrary{shadows.blur}
\usetikzlibrary{positioning}
\usetikzlibrary{backgrounds}
\usetikzlibrary{matrix}

\definecolor{apply}{rgb}{0.3,.7,0.}
\definecolor{applyface}{rgb}{.9,.95,0.}
\definecolor{invert}{rgb} {0.5,0.,1.}
\colorlet{invertface}{invert!60!red}

\usepackage{pgfplots}
\usepgfplotslibrary{groupplots}

\newsiamthm{claim}{Claim}
\newsiamremark{remark}{Remark}
\newsiamremark{hypothesis}{Hypothesis}
\crefname{hypothesis}{Hypothesis}{Hypotheses}

\usepackage{algpseudocodex}

\usepackage{graphicx,epstopdf}

\Crefname{ALC@unique}{Line}{Lines}

\usepackage{amsopn}

\usepackage{xspace}
\usepackage{bold-extra}
\usepackage[most]{tcolorbox}

\colorlet{texcscolor}{blue!50!black}
\colorlet{texemcolor}{red!70!black}
\colorlet{texpreamble}{red!70!black}
\colorlet{codebackground}{black!25!white!25}

\newcolumntype{C}[1]{>{\centering\arraybackslash}m{#1}}

\newcommand{\mesh}{\mathbb M}

\title{An implementation of tensor product patch smoothers on GPU
\thanks{Submitted to the editors February 28, 2024.
\funding{Cu Cui was supported by the China Scholarship Council (CSC) under grant no. 202106380059}}}

\author{Cu Cui\thanks{Interdisciplinary Center for Scientific Computing (IWR), Heidelberg University, Im Neuenheimer Feld 205, 69120 Heidelberg, Germany (\email{cu.cui@iwr.uni-heidelberg.de}).}
\and Paul Grosse-Bley\thanks{Institute of Computer Engineering (ZITI), Im Neuenheimer Feld 368, 69120 Heidelberg, Germany}
\and Guido Kanschat\footnotemark[2] \and Robert Strzodka\footnotemark[3] }

\headers{An implementation of tensor product patch smoothers on GPU}{Cu Cui, Paul Grosse-Bley, Guido Kanschat and Robert Strzodka}

\ifpdf
\hypersetup{ pdftitle={An implementation of tensor product patch smoothers on GPU} }
\fi

\begin{document}
\maketitle

\begin{abstract}
We present a GPU implementation of vertex-patch smoothers for higher order finite element methods in two and three dimensions.
Analysis shows that they are not memory bound with respect to GPU DRAM, but with respect to on-chip scratchpad memory. 
Multigrid operations are optimized through localization and reorganized local operations in on-chip memory, achieving minimal global data transfer and a conflict free memory access pattern. Performance tests demonstrate that the optimized kernel is at least 2 times faster than the straightforward implementation for the Poisson problem, across various polynomial degrees in 2D and 3D, achieving up to 36\% of the peak performance in both single and double precision on Nvidia A100 GPU.
\end{abstract}



\section{Introduction}
\label{sec:intro}
Our goal in this article is the proof that the multigrid V-cycle with vertex patch smoothers, which has been optimized for numerical performance in~\cite{witte2021fast} can be implemented on current GPU accelerators and utilize a reasonable share of their computing power. To this end, we take the Poisson equation on the unit cube as the fruitfly of boundary value problems and discretize with finite elements of various order. Solution methods for this problem have been optimized in every aspect from a numerical point of view and multigrid methods with block Gauss-Seidel smoothers perform very well. But, the mathematically most efficient solvers are rarely adapted to highly parallel computing architectures. We do this here and obtain a fast solver on GPU accelerators which achieves up to 36\% of the peak performance of Nvidia A100 CUDA Cores.

Eight out of the top ten supercomputers (according to the HPL benchmark) employ many-core accelerators, seven of which are GPUs~\cite{top500June23}. With the more communication dominated HPCG benchmark this number increases to nine out of the top ten \cite{hpcgJune23}. 
When ranked according to energy efficiency (with HPL again) all of the top ten most energy efficient supercomputers are GPU based \cite{green500June23}. Therefore, it is imperative to invest in the development of algorithms specifically tailored for these throughput-oriented architectures. This investment is critical to ensure efficient utilization of computational resources not only now, but also in the future.

It is well known that multigrid methods are very efficient for solving large linear systems and exhibit asymptotic complexity $\mathcal{O}(N)$, where $N$ is the number of degrees of freedom (DoF) in the system. When applied to higher order finite element methods, the choice of smoothing operator is of paramount importance in order to realize an efficient multigrid method.  \cite{kanschat2008robust} indicates that a block Gauss-Seidel smoother provides an efficient and robust multilevel preconditioner for advection–diffusion problems. Additionally, for solving Stokes equations, an overlapping Schwarz-type smoother has been demonstrated to be a valuable choice~\cite{kanschat2015multigrid}. Although these methods are highly efficient numerically requiring only few iteration steps, the computational cost in terms of time and memory can become impractical when standard computations of inverse matrices are employed to solve local problems. In~\cite{witte2021fast}, fast diagonalization was applied to the local solvers, reducing the number of operations considerably. Thus, there exists a numerically very efficient solver, which exhibits a complex algorithmic structure.

As the problem size expands, matrix-based methods ---even with sparse storage schemes--- become less efficient due to slow access to main memory.
In practice, iterative solvers often do not require the explicit form of the matrix, but rather its product with vectors, which allows fast evaluation of the discrete operator on-the-fly.
Employing sum-factorization techniques, matrix-free methods have demonstrated a significant improvement in performance, approaching near peak performance on GPU~\cite{swirydowicz2019acceleration}, compared to the implementation based on sparse matrix formats and have been successfully applied to a wide range of finite element operators \cite{KronbichlerKormann12, KronbichlerKormann19, muthing2017high, pazner2020efficient}.
Matrix-free methods have been studied extensively on CPUs and  are readily accessible through numerous open-source finite element libraries like deal.II~\cite{dealII94}, DUNE~\cite{bastian2021dune} and MFEM~\cite{anderson2021mfem}.

Implementation of matrix-free methods on GPUs have focused more on the evaluation of finite element operators, such as hyperbolic systems of conservation laws~\cite{klockner2009nodal}, acoustic and isotropic elastic wave models~\cite{modave2016gpu}.
So far, geometric multigrid methods on GPU have been implemented for finite volume methods~\cite{pazner2020efficient}, finite difference methods~\cite{feng2014numerical}. Point smoothers for finite elements were studied on GPUs in~\cite{KronbichlerLjunqkvist19}.
Nvidia's AmgX library~\cite{naumov2015amgx} provides algebraic multigrid methods, yet it relies on stored sparse matrices. 
Additionally, software tools like ExaStencils~\cite{lengauer2020exastencils} and PETSc~\cite{mills2021toward} support GPU but are constrained to low-order stencils.

Following the ideas in article~\cite{witte2021fast}, we introduce an implementation of tensor product patch smoothers specifically tailored to GPUs.
Similar to~\cite{Wichrowski2023}, we have localized and fused the operations within the smoothing operator, minimizing global vector access and eliminating the need for device synchronization.
Additionally, we conduct a thorough analysis of on-chip memory access, delving into the effects of memory and computation layouts on performance, and put forth a conflict-free memory access pattern.
These optimizations, tailored to the GPU architecture, in conjunction with mathematical refinements in the local solver, render our proposed algorithms highly efficient, offering very good performance for both low-order and high-order elements.

This article is organised as follows. In~\Cref{sec:model}, we describe the model problem and outlines the multigrid setting, along with some convergence results. 
A detailed GPU implementation is then presented in~\Cref{sec:gpu_impl}. \Cref{sec:perf} evaluates the performance of our kernels and further optimisation choices. Finally, we conclude this work in~\Cref{sec:conclusion}.

\section{Model problem and multigrid setting}
\label{sec:model}

In this section, we introduce notation and the model problem we are going to solve as well as the algorithms and implementation details which are not specific to implementation on the GPU.
We consider finite element discretizations of Poisson's equation
\begin{gather}\label{eq:model}
\begin{aligned}
    -\Delta u &= f & \text{ in } & \Omega \\
    u &{} = 0 & \text{ on } &\partial\Omega
\end{aligned}
\end{gather}
with Dirichlet boundary conditions on a rectangular domain $\Omega \subset \mathbb{R}^d$ with dimension $d = 2,3$. While this is one of the simplest and most studied problems in partial differential equations, we choose it particularly because highly optimized algorithms for its solution exist.
Multigrid methods provide for convergence independent of mesh size~\cite{Bramble93,Hackbusch85} and vertex patch smoothers exhibit fast convergence for any polynomial degree. An efficient implementation with fast diagonalization~\cite{witte2021fast} is available. On the other hand, fast implementations for multicore architectures~\cite{KronbichlerKormann12} and for GPU architectures~\cite{KronbichlerLjunqkvist19} for some of the subproblems are available. Hence, our goal is to combine those features and compete in a favorable way. 

The domain is covered by a sequence of Cartesian meshes,
\begin{gather*}
    \mesh_0 \sqsubset \mesh_1\dots\sqsubset
    \mesh_\ell\sqsubset\dots\sqsubset\mesh_L,
\end{gather*}
where $\mesh_0$ consists of a single rectangular cell, namely $\Omega$ itself. The other meshes are defined recursively by refinement. To this end, the edges of a cell of $\mesh_\ell$ are cut in half and the divided edges are reconnected to $2^d$ cells of $\mesh_{\ell+1}$.

Next, we define a finite element space $V_\ell$ on $\mesh_\ell$.
We employ a conforming finite element method with tensor product shape functions of degree $k$ denoted by $\mathbb{Q}_k$ to approximate solutions to this problem. Our method is not restricted to Lagrangian finite elements, but it is easiest explained for them. To this end, introduce a set of $k+1$ support points on the reference interval $[0,1]$, typically Gauss-Lobatto points for numerical stability. On a mesh cell $K\in\mesh_\ell$, we define support points by mapping the reference set to the dimensions of the cell and taking the $d$-fold tensor product, obtaining $(k+1)^d$ support points. A basis for $\mathbb Q_k$ on the cell $K$ is obtained by Lagrange interpolation with respect to these support points.

A basis $\{\phi_i\}_{i=0,\dots,k}$ for the finite element space $V_\ell$ is usually, see for instance~\cite{ciarlet2002finite}, obtained from local bases on each mesh cell $K$.
Note that the basis functions associated with support points in the interior of $K$ have support in $K$ and vanish outside of $K$. Hence, they can also serve as basis functions for $V_\ell$. On the other hand, basis functions associated with support points on the boundary of $K$ must be concatenated with such functions on the neighboring cells to generate continuous basis functions of $V_\ell$ which have support in several neighboring cells. Basis functions associated with the support points at the boundary of the domain are eliminated, hence implementing the boundary condition $u=0$.

The discretized weak formulation of equation~\cref{eq:model} reads as follows: Find $u_\ell \in V_\ell$ such that
\begin{equation}
  \label{eq:weak_form}
     \int_{\Omega} \nabla u_\ell \cdot \nabla v_\ell
     \,d\mathbf{x}
     =\int_{\Omega} f v_\ell
     \,d\mathbf{x}
     \quad \forall v_\ell \in V_\ell.
\end{equation}
Entering the basis functions $\{\phi_i\}_{i=0,\dots,k}$ of $V_\ell$
and the representation $u_\ell =\sum x_i \phi_i$ with the coefficient vector $x_\ell$
into~\cref{eq:weak_form}, we obtain the linear system 
\begin{equation}
    \label{eq:linear_system}
    A_{\ell} x_{\ell} = b_{\ell}.
\end{equation}
 The entries $a_{ij}$ of the matrix $A_\ell$ and $b_i$ of the right hand side vector $b_\ell$ are defined as
\begin{gather}
\label{eq:matrix-def}
    a_{ij} = \int_{\Omega} \nabla \phi_j \cdot \nabla \phi_i
     \,d\mathbf{x},
     \qquad
     b_i = \int_{\Omega} f \phi_i \,d\mathbf{x}.
\end{gather}

\subsection{Geometric Multigrid Methods}
\label{sec:multigrid}

Since the finite element spaces $V_\ell$ are nested, we define the prolongation operator $I_{\ell}^{\uparrow}$ as the embedding operator from $V_{\ell}$ to $V_{\ell+1}$ implemented as polynomial interpolation from coarser to finer mesh. The restriction operator $I_{\ell}^{\downarrow}\colon V_{\ell+1}\to V_\ell$ is defined as the transpose of the prolongation operator.

The V-cycle, see for instance~\cite{Bramble93,BrandtLivne11,Hackbusch85}, as the main component of the multigrid algorithm is described in~\cref{alg:v_cycle}.
\begin{algorithm}[tp]
\caption{Multigrid V-cycle on level $\ell$.}\label{alg:v_cycle}
\begin{algorithmic}[1]
\Procedure{${x_\ell = \textsc{Vcycle}_\ell}$}{$A_\ell,x_\ell,b_\ell$}
\If{$\ell = 0$}
    \State \Return $x_0 \gets A_0^{-1}b_0$ \Comment{coarse grid solver}
\EndIf
\State $x_\ell \gets S_\ell(x_\ell,b_\ell)$ \Comment{pre-smoothing}
\State $r_\ell \gets b_\ell - A_\ell x_\ell$ \Comment{residual}
\State $r_{\ell-1} \gets I_{\ell-1}^\downarrow r_\ell$ \Comment{coarsen} 
\State $x_\ell \gets x_\ell
    + I_{\ell-1}^\uparrow \textsc{Vcycle}_{\ell-1}(A_{\ell-1},0,r_{\ell-1})$ \Comment{coarse grid correction}
\State \Return $x_\ell \gets S_\ell(x_\ell,b_\ell)$ \Comment{post-smoothing}
\EndProcedure
\end{algorithmic}
\end{algorithm}
In addition to the level matrices $A_\ell$ and the prolongation and restriction operators, it employs the smoothing operator $S_l$, which we will discuss it in detail in the next section. In order to simplify the discussion, we will always assume a single pre- and post-smoothing step. The V-cycle can be used either as a preconditioner in a Krylov space method or as an iterative solver by itself. In the latter case, we apply it as component in the full multigrid solver~\cref{alg:linear_MG}.

\begin{algorithm}[tp]
\caption{Full multigrid solver.}\label{alg:linear_MG}
\begin{algorithmic}[1]
\State Solve $A_0x_0 = b_0$ \Comment{coarse grid solver}
\For{$\ell=1,\dots,L$}
\State $x_\ell \gets I^\ell_{\ell-1} x_{\ell-1}$ \Comment{prolongate}
\State $x_\ell \gets \textsc{Vcycle}_\ell(A_\ell, x_\ell,b_\ell)$ \Comment{V-cycle on level $\ell$}
\EndFor
\State $\delta \gets \delta_0 \coloneqq \|b_L\|$ \Comment{initial norm}
\While{$\delta > \epsilon\delta_0$}
\State $x_L \gets \textbf{V}_L(A_L,x_L,b_L)$ \Comment{V-cycle on the highest level $L$}
\State $\delta \gets \|b_L - A_Lx_L\|$ \Comment{recompute norm}
\EndWhile
\end{algorithmic}
\end{algorithm}
While the V-cycle starts on the finest level, the full multigrid solver starts on the coarsest level and successively refines the solution, before it iteratively improves it by the V-cycle until a prescribed tolerance is reached. The multigrid recursion must be closed by a coarse-grid solver on $\mesh_0$; here, we simply use one step of the smoother on level 0, where it is an exact solver.
When we report iteration steps for the full multigrid method, the number refers to the repetition of the while loop at line 6 in~\cref{alg:linear_MG}.

\subsection{Vertex-patch smoothers}
\label{sec:vps}

Let $\mathbb V_\ell$ be the set of all interior vertices of the mesh $\mesh_\ell$ on level $\ell$ of the multilevel discretization. For any $p_j \in \mathbb V_\ell$, the corresponding vertex patch consists of all cell having $p_j$ as one of their vertices. For regular quadrilateral and hexahedral meshes, the vertex patches consist of 4 and 8 cells, respectively. The set of all vertex patches forms an overlapping covering of the whole mesh $\mesh_\ell$. For simplicity of notation, we will drop the level index $\ell$ in the discussion of the smoother and introduce $J$ as the cardinality of $\mathbb V_\ell$.

\begin{figure}[tp]
\centering
    \includegraphics[width=.45\textwidth]{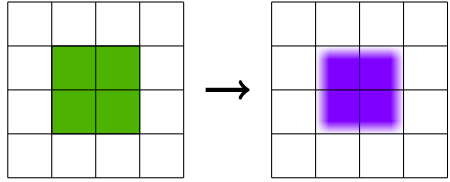}
  \caption{Degrees of freedom including the boundary of the vertex patch, defining $\overline V_j$, and used for the local computations of residuals (left). Interior degrees of freedom involved in the local solver, defining $V_j$ (right).}
  \label{fig:vertex_patches}
\end{figure}

The vertex-patch smoother is a domain decomposition method which uses all vertex-patches as its subdomains. Therefore, we associate with each vertex patch $j$ the subspace $V_j$ of the finite element solution space consisting of functions with support on the patch and being zero elsewhere. The local correction of a current approximate solution $x$ then consists of two steps, first computing the residual
\begin{gather}
    \label{eq:residual}
    r \gets b-Ax,
\end{gather}
and the local solver
\begin{gather}
\label{eq:local-solver}
    y \gets x+R_j^T A_j^{-1}R_j r.
\end{gather}
Here, the restriction operator $R_j$ selects only the coefficients from the vector $r$ which are associated with basis functions in $V_j$. Similarly, the matrix $A_j$ is obtained by restricting the definition of the matrix in~\cref{eq:matrix-def} to the basis functions spanning $V_j$.
We call the mapping from $x$ to $y$ introduced by equations~\cref{eq:residual} and~\cref{eq:local-solver} the subspace correction $y = s_j(x)$ and define the multiplicative vertex-patch smoother by
\begin{gather}
    \label{eq:simple-loop}
    x \gets s_j\circ s_{j-1} \circ\dots\circ s_1(x)
    .
\end{gather}
From the point of view of implementation, this algorithm has two drawbacks. First, every step involves a multiplication with the global matrix $A$ to compute the residual, thus resulting in a complexity quadratic in the number of mesh cells. Second, due to the overlap of the vertex-patches, it is strictly sequential. We address the first problem by observing that the operation $R_j Ax$, which combines parts of~\cref{eq:residual} and~\cref{eq:local-solver}, involves only those entries of $x$ associated to the basis functions of $V_j$ and those of basis functions whose support overlaps with the vertex-patch $j$. These correspond to the subspace $\overline V_j$, which is spanned by the basis functions associated with all degrees of freedom of the vertex-patch, including those on the boundary. This situation is visualized in~\cref{fig:vertex_patches}. 
Hence, we obtain an equivalent definition of $s_j$ by
\begin{gather}
    r_j \gets R_j b - \overline A_j \overline R_j x,
    \qquad
    y \gets x + R_j^T A_j^{-1} r_j,
\end{gather}
where $\overline A_j$ is the rectangular matrix obtained by constraining $\phi_i$ and $\phi_k$ in~\cref{eq:matrix-def} to the basis functions of the spaces $V_j$ and $\overline V_j$, respectively. This version of $s_j$ only involves degrees of freedom local to the vertex-patch. The impact of this localization on the performance of the vertex-patch smoother on a CPU is studied in~\cite{Wichrowski2023}.

\begin{figure}[tp]
\centering
    \includegraphics[width=.65\textwidth]{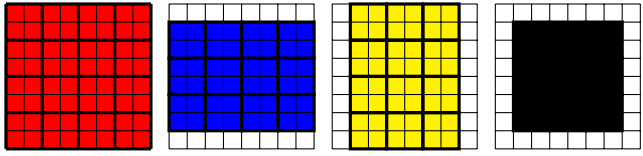}
    \caption{Non-overlapping coloring for vertex patches on regular meshes. In each color we strive to obtain a parqueting of the whole domain. Subsequent colors are obtained by shifting the first patch by one cell in each coordinate direction.}
    \label{fig:patch_coloring}
\end{figure}

The localization of the residual computation just discussed avoids unnecessary computations. But it has a second consequence. Referring to~\cref{fig:vertex_patches} again, it is obvious that the update produced by one patch affects the update of another only if they are overlapping. Hence, we can use this property to subdivide the set of patches. Indeed, all non-overlapping patches in~\cref{fig:patch_coloring} can be processed in parallel without causing any races.
We follow the custom in the literature and call this subdivision of the patches ``colorization''. Summarizing this section, we present the colorized vertex-patch smoother in ~\cref{alg:Multiplicative_Schwarz}.
\begin{algorithm}[tp]
\caption{Colorized Vertex-Patch Smoother $S_\ell(x_\ell,b_\ell)$.}
\label{alg:Multiplicative_Schwarz}
\begin{algorithmic}[1]
\For{$c=1,\dots,n_{c}$} \Comment{Sequential loop over colors}
\For{\textbf{each } $j \in\mathcal{J}_c$} \Comment{parallel loop inside color}
\State $r_j \gets R_j b - \overline A_j \overline R_j x_j$ \Comment{local residual}
\State $x_\ell \gets x_\ell + R^T_j A_j^{-1} r_j$ \Comment{local solver}
\EndFor
\EndFor
\end{algorithmic}
\end{algorithm}
There, the set of patches has been split up in $n_c$ subsets represented by the index sets $\mathcal J_c \subset \{1,\dots,J\}$ where $c=1,\dots,n_c$. The simple loop of equation~\cref{eq:simple-loop} has been replaced by two nested loops. The outer loop is sequential and runs over the different ``colors''. The inner loop is parallel and runs over all vertex-patches of a single color.

The vertex-patch smoother is fairly complex and expensive. Note that for tensor elements $\mathbb Q_k$, a linear system of dimension $(2k-1)^d$ must be solved for every vertex patch~\cite{braess2007finite}, which is prohibitively expensive if standard techniques like QR factorization are used. Before improving it algorithmically, let us justify its use by a few convergence results.

 To justify its use, we present results in~\cref{tab:LMG_v_2d,tab:LMG_v_3d}. We solve the linear system~\cref{eq:linear_system} in two and three dimensions for right hand side $f \equiv 1$ and for varying finest level $L$ to a relative accuracy of $10^{-9}$. As the results show, the number of iteration steps is independent of the mesh level except for lowest order, where it saturates more slowly in two dimensions. Moreover, from degree 3 onwards, at most 3 steps of the method are needed to convergence, making it almost a direct solver.
\begin{table}[tp]
    \centering
    \small
    \caption{Iteration steps $n$ for full multigrid with vertex-patch smoother in two dimensions for different mesh levels $L$ and polynomial degrees $k$. Relative accuracy gain of $10^{-9}$.}
    \label{tab:LMG_v_2d}
    \begin{tabular}%
    {c|cccccccccc}
        \hline
         \textbf{$L$} & $\mathbb{Q}_1$ & $\mathbb{Q}_2$ & $\mathbb{Q}_3$ &
         $\mathbb{Q}_4$ & $\mathbb{Q}_5$ & $\mathbb{Q}_6$ & $\mathbb{Q}_7$ &
         $\mathbb{Q}_8$ & $\mathbb{Q}_9$ & $\mathbb{Q}_{10}$ \\
        \hline
         4 & 9 & 5 & 3 & 3 & 3 & 2 & 2 & 2 & 2 & 2 \\[-1.5ex]
    \vdots & \vdots & \vdots & \vdots & \vdots & \vdots & \vdots & \vdots & \vdots & \vdots & \vdots \\[-1.5ex]
         11 & 7 & 5 & 3 & 3 & 3 & 2 & 2 & 2 & 2 & 2 \\
         12 & 7 & 4 & 3 & 3 & 2 & 2 & 2 & 2 & 2 & 2 \\
        \hline
    \end{tabular}
\end{table}
\begin{table}[tp]
    \centering
    \small
    \caption{Iteration steps $n$ for full multigrid with vertex-patch smoother in three dimensions for different mesh levels $L$ and polynomial degrees $k$. Relative accuracy gain of $10^{-9}$. Dashes are not computed due to excessive memory.}
    \label{tab:LMG_v_3d}
    \begin{tabular}%
    {c|cccccccc}
        \hline
         \textbf{$L$} & $\mathbb{Q}_1$ & $\mathbb{Q}_2$ & $\mathbb{Q}_3$ &
         $\mathbb{Q}_4$ & $\mathbb{Q}_5$ & $\mathbb{Q}_6$ & $\mathbb{Q}_7$ &
         $\mathbb{Q}_8$ \\
         \hline
         4  &  6  &  5 &  3  & 3  & 3  & 3 & 2 & 2 \\[-1.5ex]
         \vdots & \vdots & \vdots & \vdots & \vdots & \vdots & \vdots & \vdots & \vdots \\[-1.5ex]
         7  &  6  &  5 &  3  & 3  & 3  & 3 & 2 & 2  \\
         8  &  6  &  5 &  3  & 3  & 3  & 3 & 2 & ---  \\
         9  &  6  &  5 &  3  & ---  & ---  & --- & --- & --- \\
         10 &  6  &  --- &  ---  & ---  & ---  & --- & --- & --- \\
        \hline
    \end{tabular}
\end{table}

In~\cref{tab:vp_GS_LMG}, we compare iteration steps of the multigrid iteration with vertex-patch and point Gauss-Seidel smoothers, respectively, for higher order elements. All results were obtained on fine meshes such that the iteration steps already remained constant. We again observe that higher order elements require fewer iterations of the multigrid V-cycle compared to lower order elements. The number of iterations required for point Gauss-Seidel grows with polynomial degree and is 15x more than that for vertex-patch smoother at $\mathbb{Q}_8$.
Thus, we conclude that the multiplicative vertex-patch smoother yields a mathematically very efficient multigrid method and is essential for higher order elements.
But it is computationally expensive compared to other smoothers. Therefore, efficient implementation is needed to make it feasible.
\begin{table}[tp]
    \centering
    \small
    \caption{Comparison of vertex patch smoother and point Gauss–Seidel smoother for higher order elements in three dimensions (iteration steps for residual reduction of $10^{-9}$ on mesh levels, where they remain constant).}
    \label{tab:vp_GS_LMG}
    \begin{tabular}{lC{0.7cm}C{0.7cm}C{0.7cm}C{0.7cm}C{0.7cm}C{0.7cm}C{0.7cm}C{0.7cm}}
        \hline
           & $\mathbb{Q}_1$ & $\mathbb{Q}_2$ & $\mathbb{Q}_3$ &
         $\mathbb{Q}_4$ & $\mathbb{Q}_5$ & $\mathbb{Q}_6$ & $\mathbb{Q}_7$ &
         $\mathbb{Q}_8$ \\
         \hline
         point GS & 6 & 8 & 11 & 13 & 19 & 20 & 27 & 30 \\
         vertex patch & 6 & 5 & 3 & 3 & 3 & 3 & 2 & 2 \\
        \hline
    \end{tabular}
\end{table}

\subsection{Matrix-Free Operation}
\label{sec:mf}
For finite element codes on modern hardware, the transfer of data from memory to the processing unit is typically the bottleneck limiting processing speed, this is known as the memory wall. This is true for implementations on CPU and GPU, see for instance~\cite{KronbichlerKormann12,KronbichlerKormann19,KronbichlerLjunqkvist19}. Thus, we decide on a matrix free implementation~\cite{dealII94}, where the matrices $A_\ell$ are never stored, but reconstructed locally on a cell from the bilinear form on the left of~\cref{eq:weak_form} whenever needed.

Considering a matrix-free implementation for computing the local stiffness matrix on a two-dimensional cell K, we have:
\begin{equation}\label{eq:mf0}
    A_{ij}^K = \sum_q \hat\nabla \hat\phi_i^T(\hat x_q) J^{-T}(\hat x_q) |J(\hat x_q)|
     J^{-1}(\hat x_q) \hat\nabla\phi_j(\hat x_q) w_q,
\end{equation}
where $\hat\nabla$ is the gradient on the reference cell, $J(\hat x_q)$ the Jacobian of the transformation from the reference to the real cell at point $\hat x_q$ and $w_q$ the quadrature weight. The inverse Jacobian on a Cartesian mesh is both diagonal and constant, eliminating the need to store individual Jacobians for each quadrature point on every cell and simplifying storage. With a basis constructed from the tensor product of one-dimensional shape functions, sum factorization can then be used to evaluate all basis functions at all of the quadrature points efficiently~\cite{orszag1979spectral}. Therefore, the local operator~\cref{eq:mf0} is equivalent to
\begin{equation}\label{eq:mf}
    A^K = \mathbf{G}^T \mathbf{W} \mathbf{G},
\end{equation}
where $\mathbf{G}=\left[\begin{array}{c}S \otimes D \\ D \otimes S \end{array}\right]$, $S$ is the matrix of all one-dimensional shape functions evaluated at quadrature points and $D$ the matrix of their derivatives and $W_{i_q,j_q} = J^{-T}(\hat x_{i_q,j_q}) |J(\hat x_{i_q,j_q})| J^{-1}(\hat x_{i_q,j_q}) w_{i_q,j_q}$. As the present work is focused on Cartesian geometry, we consider a specific optimization $A^K = A_1 \otimes M_0 + M_1 \otimes A_0$, which is used in the evaluation of smoothing operation.

\subsection{Fast Diagonalization}
\label{sec:fd}
With the use of tensor product elements, the local discretization $\overline A_j$ in~\cref{alg:Multiplicative_Schwarz} on a Cartesian mesh has a separable Kronecker representation of the form:
\begin{equation}
\label{eq:tensorproduct_2d}
    \overline A_j = A_j^1 \otimes M_j^0 + M_j^1 \otimes A_j^0
\end{equation}
in two dimensions and
\begin{equation}
\label{eq:tensorproduct_3d}
    \overline A_j = A_j^2 \otimes M_j^1 \otimes M_j^0 + M_j^2 \otimes A_j^1 \otimes M_j^0 + M_j^2 \otimes M_j^1 \otimes A_j^0
\end{equation}
in three dimensions, where $A_j^d$ and $M_j^d$ are one dimensional stiffness and mass matrices on patch $j$, respectively.
Since $A_j$ has an identical structure just with smaller matrices, its exact inverse can be found through the fast diagonalization method~\cite{lynch1964direct}:
\begin{equation}
\label{eq:inverse2d}
    A_j^{-1} = \bigl[S_j^1 \otimes S_j^0\bigr]
    \bigl[\Lambda_j^1 \otimes I + I \otimes \Lambda_j^0\bigr]^{-1} 
    \bigl[S_j^1 \otimes S_j^0\bigr]^T
\end{equation}
in two dimensions and
\begin{equation}
\label{eq:inverse3d}
    A_j^{-1} = \bigl[S_j^2 \otimes S_j^1 \otimes S_j^0\bigr]
    \bigl[\Lambda_j^2 \otimes I \otimes I + I \otimes \Lambda_j^1 \otimes I + I \otimes I \otimes \Lambda_j^0\bigr]^{-1}
    \bigl[S_j^2 \otimes S_j^1 \otimes S_j^0\bigr]^T
\end{equation}
in three dimensions, where $S_d$ is the matrix of eigenvectors to the generalized eigenvalue problem in the given tensor direction $d$:
\begin{equation}
\label{eq:fast_inverse}
    A_j^d S_j^d=M_j^d S_j^d\Lambda_j^d , \quad d=0,...,\mathrm{dim},
\end{equation}
and $\Lambda_d$ is the diagonal matrix representing the generalized eigenvalues $\lambda_i$.
By exploiting the tensor structure elements, the local solver $A_j^{-1}$ is performed at the cost of $\mathcal{O}(dk^{d+1})$, which is significantly reduced compared to explicit inversion $\mathcal{O}(k^{3d})$.
This is the same asymptotic cost as sum factorization.
Another benefit is the significant decrease in memory consumption, this algorithm requires to store only one dimensional eigenvectors and eigenvalues. Numerical experiments conducted in the next section demonstrate that this is the key to achieving an efficient smoother.
\section{GPU Implementation}
\label{sec:gpu_impl}

In this section, we explore the GPU implementation of the multigrid method outlined in~\cref{alg:v_cycle}. It is based on the Compute Unified Device Architecture (CUDA) programming model provided by NVIDIA for developing general purpose computing kernels that run on their massively parallel GPUs~\cite{cuda2023}. 

The multigrid V-cycle involves several core operations. Instead of using a single large kernel that includes all operations, we divide the different operations into several kernels when designing the GPU algorithm. This makes it easier to identify bottlenecks in the algorithm and provides flexibility to optimize each kernel based on its own characteristics.

We consider the following building blocks for implementation:
the {\it operator application} which performs the global matrix-vector product as defined in~\cref{eq:linear_system}; the {\it grid transfer operation} which performs the prolongation and restriction between levels; the {\it smoothing operation} which performs the Multiplicative Schwarz Smoother according to~\cref{alg:Multiplicative_Schwarz} and the {\it vector operations} which compute element-wise and reduction-type vector operations. These operations are launched sequentially by the CPU, allowing for easy integration with existing CPU frameworks. The massively parallel acceleration is achieved within each operation, leveraging the GPU's processing power efficiently. By structuring the implementation this way, the migration from CPU code to GPUs becomes more manageable and allows the focus to be on optimizing and parallelizing the individual kernels for better performance. 

The first three operations in the multigrid method are typically implemented by traversing all cells (or patches in the vertex-patch smoother). The operations performed on each cell are independent of others and can be executed in parallel. Within each cell, we evaluate the values of all degrees of freedom simultaneously without any communication with other cells.
This computation pattern is well-suited for the CUDA programming model. Each CUDA thread is associated with a degree of freedom, and each CUDA thread block with a cell, thus allowing for two levels of parallelism which enables more efficient computation. When designing the CUDA kernel, the focus is on the thread block. In ~\cref{alg:cuda_pattern} we describe a typical pattern for CUDA programming. A key to achieving high performance on the GPU is to make good use of the limited on-chip memory~\cite{kirk201771}. To this end, we analyze the implementation with focus on on-chip memory\footnote{On-chip memory in this work refers to \textit{Shared Memory} in CUDA programming~\cite{cudaProg}. It is a low latency on-chip memory like L1-cache that is visible to all threads in a thread block and is much faster than global memory (DRAM).} usage. This ensures that the most critical data is stored on-chip, reducing expensive access to DRAM and maximizing computation efficiency.
\begin{algorithm}[tp]
\caption{A typical pattern for CUDA programming.}\label{alg:cuda_pattern}
\begin{algorithmic}[1]
\State Load data from DRAM to on-chip memory
\State Synchronize thread block to make sure the data has been loaded
\State Perform computation on data in on-chip memory
\State Synchronize thread block again to make sure the result has been updated
\State Store the result from on-chip memory back to DRAM
\end{algorithmic}
\end{algorithm}

All numerical experiments in this work are performed on a single NVIDIA Ampere A100 SXM4 GPU with 80GB of high-speed HBM2e memory for DRAM which provides 2TB/s peak memory bandwidth, hosted on a system with two AMD EPYC 7282 16-Core processors. 
For accurate data collection, the test is conducted 100 times, averaging every 10 runs and selecting the best time among those ten averages. This approach eliminates the impact of a cold cache and potential processor frequency changes~\cite{hoefler2015scientific}.


\subsection{Smoothing operation}
\label{sec:smoothing}
Even with only one pre- and one post-smoothing step, the smoothing with overlapping patches is still the most expensive component of the multigrid algorithm. 
Early testing indicated that performing the smoothing on the GPU and leaving the CPU for other operations was not a good option. This is because data transfer between GPU and CPU takes a long time, affecting the overall performance of the algorithm. Therefore we have implemented the entire multigrid cycle on the GPU. During the setup stage, vectors and relevant data are allocated and prepared on the GPU, and there is no data communication between the CPU and GPU throughout solving the linear system.

\subsubsection{Implementation of local operations}

In~\cref{alg:Multiplicative_Schwarz}, the local solver essentially performs a matrix-vector multiplication of precomputed dense matrix $A^{-1}$ and local vector $r$. It is clear that the matrix-vector multiplication is memory-bound on GPU and its performance is highly influenced by DRAM access patterns and data transfer sizes. By exploiting the tensor product structure, the exact inverse can be computed through fast diagonalization~\cref{eq:inverse2d,eq:inverse3d}. This approach allows for the evaluation of the local solver to be transformed into a matrix-matrix multiplication of smaller dimensions. Consequently, the computational complexity is reduced from $\mathcal{O}(k^{2d})$ to $\mathcal{O}(dk^{d+1})$~\cite{witte2021fast}. With growing polynomial degree $k$, it becomes a compute-bound algorithm now, which offers us more opportunities to increase performance. 

\begin{figure}[tp]
\centering
    \includegraphics[width=.75\textwidth]{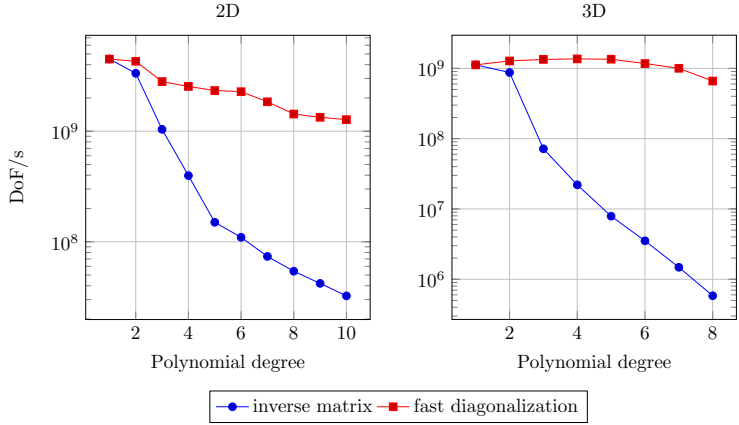}
\caption{Comparison of implementation of local operation on throughput of one smoothing step in two and three dimensions. Inverse matrix: apply the local solver through multiplication with an inverse matrix. Fast diagonalization: apply the local solver through fast diagonalization.}
\label{fig:tp_element}
\end{figure}

The throughputs of two methods with straightforward implementation are shown in~\cref{fig:tp_element}. As the problem is constant-coefficient, the inverse on all patches is the same, making it sufficient to store only one inverse. For $\mathbb{Q}_1$ elements, both methods perform identically, which is expected since only one interior node is updated. For polynomial degree $k=2$, all data fits within the on-chip memory and the observed performance enhancement is attributed to the reduced number of operations necessitated by the fast diagonalization method.
However, as the polynomial order increases, the explicitly stored matrices no longer fit in the L1 cache. Consequently, longer latency L2 cache accesses or even DRAM accesses start to dominate the performance, resulting in a rapid drop in performance for the matrix structure method. This effect is more pronounced in three dimensions. In contrast, the fast diagonalization approach reduces memory consumption from $\mathcal{O}(k^{2d})$ to $\mathcal{O}(dk^2)$, making it possible to load the required data into fast on-chip memory. This approach leads to fairly stable performance for various element orders, as indicated in~\cref{fig:tp_element}. An analogous optimization strategy can also be used for the evaluation of local residuals, where separable Kronecker representation of local matrix $\overline A_j$ can be used to avoid a straightforward matrix-vector multiplication. 

Moreover, when evaluating the tensor products~\cref{eq:tensorproduct_3d}, we can combine the evaluations of individual components. The innermost sum in x0-direction ($M_0x$) is the same for derivatives in x1- ($M_2\otimes A_1 \otimes M_0$x) and x2-direction ($A_2\otimes M_1 \otimes M_0$) in 3D, and the temporary result ($M_1\otimes A_0x$ and $A_1\otimes M_0x$) after the first summation can be reused~\cite{KronbichlerKormann12}.

\subsubsection{Computation layout and data layout}
The matrix $\overline R_j$ ($R_j$), which describes the mapping of local and global degrees of freedom, is usually stored as a vector of indices. 
By storing only the first index of each patch, the storage of $\overline{R}_j$ can be compressed into a single index. This is possible because the optimization of compressing the storage of the matrix $\overline{R}_j$ (or $R_j$) can lead to significant reductions in device memory usage. 

\begin{figure}[tp]
\centering
    \includegraphics[width=.45\textwidth]{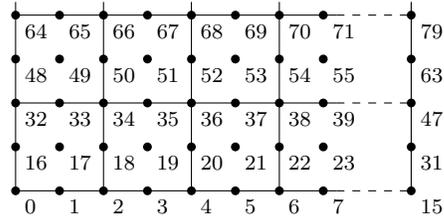}
  \caption{Visualization of degrees of freedom layout for $\mathbb{Q}_2$ element in 2D, using global lexicographical numbering.}
  \label{fig:dof_layout}
\end{figure}

This optimization is particularly applicable to a regular grid structure with global lexicographical numbering, as illustrated in ~\cref{fig:dof_layout}. Otherwise, in the case of hierarchical numbering, there is a requirement to store additional indices, encompassing one for each face as well as the interior. The use of lexicographical numbering allows for more efficient computation of each thread's index based on the first index and an appropriate stride. This index computation is considerably faster than fetching data from DRAM.
Another advantage of using lexicographical numbering is the ability to access DRAM as linearly as possible. This layout consistency facilitates more efficient memory access patterns and eliminates the need for further rearrangement during tensor product evaluation.

\begin{figure}[tp]
\centering
\includegraphics[width=.4\textwidth]{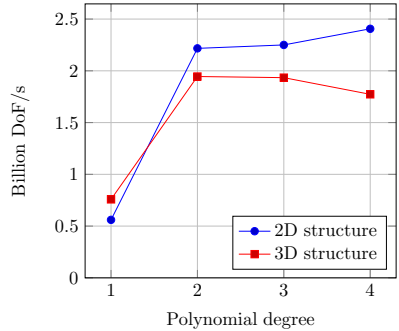}
\caption{Comparison of thread block structure of one smoothing step in three dimensions. 2D structure---each thread processes a “column” of DoFs, 3D structure---each thread processes a single DoF.}
\label{fig:2d_3d_tb}
\end{figure}

One thread per degree of freedom is a natural choice for the implementation in two dimensions. However, in three dimensions, the number of degrees of freedom in a patch can easily exceed 1024 (the maximal number of threads in a thread block), for instance, in $\mathbb{Q}_5$, we have $(2\times 5+1)^3 = 1331$ DoFs in total. Therefore, we adopt the strategy of one thread processing a `column' of $2k+1$ degrees of freedom for higher order elements.
In~\cref{fig:2d_3d_tb} we compare the one thread per degree (3D thread structure) with the one thread per `column' (2D thread structure) approach on a 3D problem. 
For higher-order elements, the 2D structure consistently outperforms. Threads with more work provide the scheduler with more opportunities to optimize, resulting in improved performance. Only in the case of polynomial degree $k = 1$ does the 3D thread structure perform better. This phenomenon arises due to the on-chip memory usage within a thread block being the constraining factor, rather than the count of threads. Consequently, increasing the thread block size serves to optimize resource utilization and enhance occupancy.

For better resource balancing among the number of registers, the number of threads, and the size of on-chip memory, the implementation processes multiple patches within each thread block. \Cref{fig:thread_layout} illustrates the thread layout and memory layout for vertex patch smoothers, where several patches are handled in one thread block. The global numbering scheme (\cref{fig:dof_layout}) we employ ensures the data is stored in a row major way. Thus, we access the memory space in a linear fashion in each row. A very important performance consideration in CUDA programming is the coalescing of DRAM accesses which minimizes the number of memory transactions needed. By employing lexicographical numbering, we ensure that DRAM accesses are coalesced whenever feasible.
It is worth noting that the data at the patch boundary is fetched twice, ensuring that several patches can be updated simultaneously without synchronization or competition. The actual number of patches per thread block depends on the precision, space dimension and element order. 
\begin{figure}[tp]
\centering
    \includegraphics[]{Figures/thread_layout.tikz}
  \caption{Thread layout and DRAM layout for vertex patch smoothers for $\mathbb{Q}_3$ element with three non-overlapping patches. $T_i$ indicates a certain memory address accessed by thread $i$. DRAM addresses are contiguous within each color.}
  \label{fig:thread_layout}
\end{figure}

\begin{figure}[tp]
\centering
\includegraphics[width=.75\textwidth]{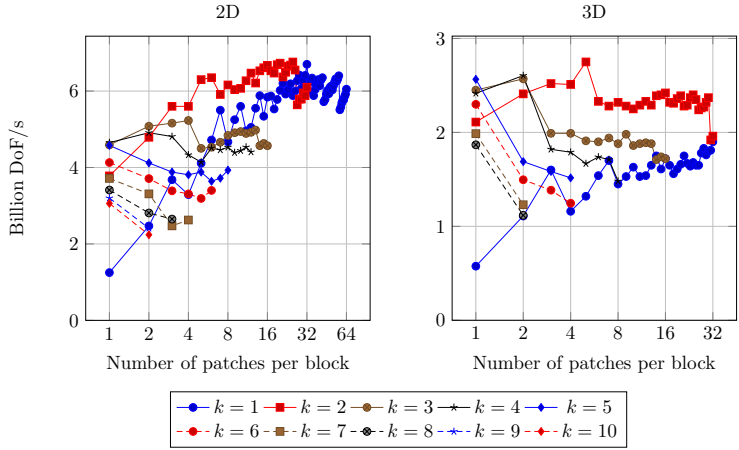}
\caption{Throughput of one smoothing step for various patches per thread-block for different polynomial degree $k$ in 2D and 3D.}
\label{fig:patches_per_shared_mem_detail}
\end{figure}

\begin{figure}[tp]
\centering
\includegraphics[width=.75\textwidth]{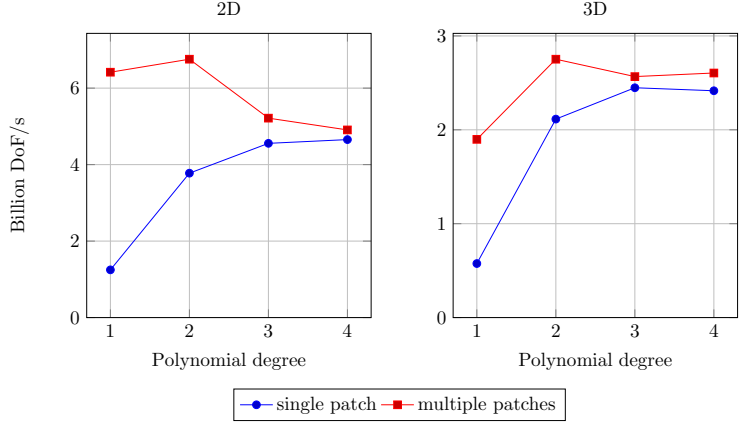}
\caption{Throughput of smoothing step for one or more patches per thread-block for lower polynomial degree in 2D and 3D.}
\label{fig:patches_per_shared_mem}
\end{figure}
Since all computations are performed in on-chip memory, the amount of on-chip memory is the decisive factor in determining the maximum number of patches each block can process.
\Cref{fig:patches_per_shared_mem_detail} displays the throughput of one smoothing step as a function of the number of patches per thread block. More patches per thread block clearly help for degrees 1 or 2, for degrees 3 or 4 there is only little benefit and higher degrees better stick to one patch per thread block. The performance results are summarized in~\Cref{fig:patches_per_shared_mem} using optimal patches per thread block for lower order elements and single-patch strategy tends to work better for polynomial degrees greater than four.

\subsubsection{Kernel Choices}

\begin{table}[tp]
    \centering
    \caption{Pseudo code of different versions of vertex-patch smoother and memory requirement. $*_\ell$ denotes the global vector on level $\ell$. $*_r$ denotes the restriction of global level vector $*_\ell$ to a specific patch $j$.}
    \footnotesize
    \label{tab:pseudo_vps}
    \begin{tabular}{|llcc|}
        \hline
         Global &  & global vector access & on-chip memory usage \\
        \hline
         residual &  $t_\ell \leftarrow A_\ell x_\ell$ & $1w + 1r$ & - \\
                        &  $r_\ell \leftarrow b_\ell - t_\ell $ & $1w + 2r$ & - \\
        \hdashline[1.5pt/1.5pt]
         solver & $r_j \leftarrow R_j r_\ell$ & $0w + 1r$ & \multirow{3}{*}{$3(2k-1)^d + 2k(2k-1)^2$}\\
                    & $t_j \leftarrow A_j^{-1}r_j $ & $0w + 0r$ &  \\
                    & $x_\ell \leftarrow x_\ell + R^T_j t_j $ & $1w + 1r$ &  \\
        \hline
        Separate &  & global vector access & on-chip memory usage \\
        \hline
         residual & $x_j \leftarrow \overline R_j x_\ell$, $b_j \leftarrow \overline R_j b_\ell$ & $0w+2r$ & \multirow{3}{*}{$(2+d-1)(2k-1)^d + 2(2k+1)^2$} \\
            & $b_j \leftarrow b_j - \overline A_j x_j$ & $0w+0r$ &  \\
                        &  $r_\ell \leftarrow \overline R^T_j b_j $ & $1w + 0r$ &  \\
        \hdashline[1.5pt/1.5pt]
         solver & $r_j \leftarrow R_jr_\ell $ & $0w+1r$ & \multirow{2}{*}{$2(2k-1)^d + 2k(2k-1)^2$} \\
                    & $x_\ell \leftarrow x_\ell + R^T_j A_j^{-1} r_j $ & $1w + 1r$ &  \\
        \hline
        Fused &  & global vector access & on-chip memory usage \\
        \hline
         residual &  $x_j \leftarrow \overline R_j x_\ell$, $b_j \leftarrow \overline R_j b_\ell$ & $0w + 2r$ & \multirow{4}{*}{$(2+d-1)(2k-1)^d + 2(2k+1)^2$} \\
                        &  $ b_j \leftarrow b_j - \overline A_j x_j $ & $0w + 0r$ &  \\
         solver & $x_j \leftarrow x_j + A_j^{-1}b_j $ & $0w+0r$ &  \\
                    & $x_\ell \leftarrow  R^T_j x_j $ & $1w + 0r$ &  \\
        \hline
    \end{tabular}
\end{table}
With all the details prepared, we can now discuss how to structure the smoother kernel. We propose three kernels including different levels of optimization. \Cref{tab:pseudo_vps} shows the pseudo code of different smoothing kernels. The {\it Global} version computes the residual globally with the matrix-free finite element operator $A_\ell$ and a separate kernel is used for the local solver. We use \texttt{LinearAlgebra::distributed::Vector} class from \texttt{deal.II} library~\cite{dealII94} for global vector operations. The data loaded into on-chip memory includes two local vectors and a temporary buffer ($3(2k-1)^d$), as well as 1D eigenvectors and eigenvalues for Fast Diagonalization ($2k(2k-1)^2$). Since $\overline A_j$ is a local operator on patches, there is no need for the global operation and its evaluation can be interleaved with vector updates, which allows further optimization. In the {\it Separate} version, the host launches two separate kernels, computing residual and inversion locally one after the other. Furthermore, local vector updates in on-chip memory can be performed in place requiring fewer temporary vectors. Because both operations are patch-based, we can easily combine them into a single kernel, as illustrated in {\it Fused} version. This avoids device synchronization while also reducing global vector access for loading and storing the residual vector.
However, the dimensions of the matrix multiplications in the two operations do not coincide, as explained in~\cref{fig:vertex_patches}. This means that some threads will not participate in the computation and will remain idle during the inversion operation. We mask the threads as inactive that correspond to the boundary degrees of freedom.

\begin{figure}[tp]
\centering
\includegraphics[width=.75\textwidth]{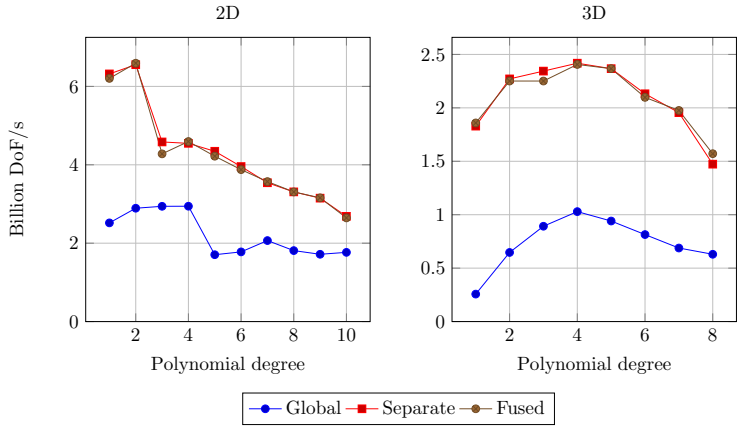}
\caption{Comparison of throughputs for one smoothing step for different smoother kernels with 26--85 million DoFs in 2D and with 135--721 million DoFs in 3D.}
\label{fig:kernels_GSF}
\end{figure}

In addition to the theoretical analysis of the algorithm, we show in~\cref{fig:kernels_GSF} the throughput for different smoother kernels with 26–85 million DoFs in 2D and with 135–721 million DoFs in 3D, where we ensure that the problem size is large enough for the algorithm to achieve near optimal performance. A clear speedup is visible by updating the residual locally. When comparing the {\it Separate} and {\it Fused} kernels, they demonstrate comparable performance; the {\it Fused} kernel involves fewer global vector accesses, while the {\it Separate} kernel necessitates less shared memory usage and minimizes thread divergence. Consequently, in the absence of conclusive evidence regarding superior performance between the two kernels, we will illustrate our optimizations using the Fused kernel as a reference. Subsequently, we will provide comparative results between the two kernels after optimizing the {\it Separate} kernel.

\subsubsection*{Matrix-Free Multiplication}
From Section~\ref{sec:mf}, the evaluation of the finite element operator acting on a vector is of the form:
\begin{equation*}
    v_\ell = A_\ell u_\ell = \sum_{K=1}^{n_{\textrm{cells}}} R_K^T A_{\ell}^K R_K u_\ell,
\end{equation*}
the result is computed as a sum of local matrix multiplication from each cell. These local evaluations are independent of one another, making it feasible to parallelize the algorithm over cells. Similar to the parallel strategy used in the patch smoother, cells are distributed into blocks and we compute the contribution from each thread block in parallel.

For each cell, we first load local vector values into on-chip memory and then perform numerical quadrature according to formula~\cref{eq:mf}. The values and gradients of the one dimensional shape functions at quadrature points are stored in constant memory, which reduces on-chip memory demand and thus increases occupancy (more simultaneous thread blocks on GPU). In sum-factorization, the evaluation of the gradient of a finite element function based on tensor-product polynomials corresponds to a matrix–matrix product. This leads us to a finer-grained parallelism where matrix multiplication can be readily parallelized in the GPU using a one thread per entry strategy.

When performing a parallel loop over cells, the summations of local vectors to global vector can potentially cause races, since DoFs on vertices or edges will be updated by several cells. A coloring approach can be used to avoid conflicts, in which only cells without sharing DoFs will be processed simultaneously. In this work, we adopt atomic operations for resolving conflicting writes, which have been demonstrated to be faster than a coloring approach~\cite{KronbichlerLjunqkvist19}.

\subsubsection*{Grid transfer operation}
A matrix-free implementation with similar ingredients as for the finite element operator evaluation is employed for grid transfer operations. A parallel loop is issued over parent cells and perform interpolation along different dimension to get the final result. Restriction operator is defined by the transpose of prolongation operator (see~\cite{KronbichlerLjunqkvist19} for details). Grid transfer operator is less interesting to us, considering its small portion of the total runtime (as shown in~\cref{fig:breakdown_GMRES_components_ratio}), compared to the smoother and it is not the main focus of this work.

\section{Performance}
\label{sec:perf}

\begin{figure}[tp]
\centering
\includegraphics[width=.75\textwidth]{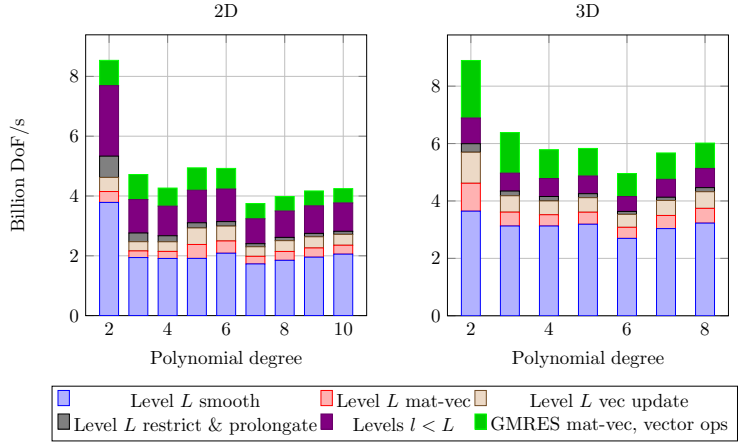}
\caption{Breakdown of components and throughput of the preconditioned GMRES solver with a vertex-patch smoother in 2D and 3D.}
\label{fig:breakdown_GMRES_components}
\end{figure}

In this section, the performance of the multigrid method is analyzed in more detail to understand the impacts of its different components in order to perform more specific optimization.
\Cref{fig:breakdown_GMRES_components} shows a breakdown of the runtime of the preconditioned GMRES solver with vertex-patch smoother in 2D and 3D with the following components:
\begin{itemize}
    \item levels $l<L$: operations on all levels except the finest one;
    \item level $L$ mat-vec: matrix-vector multiplication on the finest level $L$;
    \item level $L$ restrict \& prolongate: grid transfer operations on the finest level $L$;
    \item level $L$ smooth: pre and post-smoothing operations on the finest level $L$;
    \item level $L$ vector update: vector operations for residual on the finest level $L$;
    \item GMRES mat-vec, vector ops: all the other operations in the GMRES solver.
\end{itemize}

\begin{figure}[tp]
\centering
\includegraphics[width=.75\textwidth]{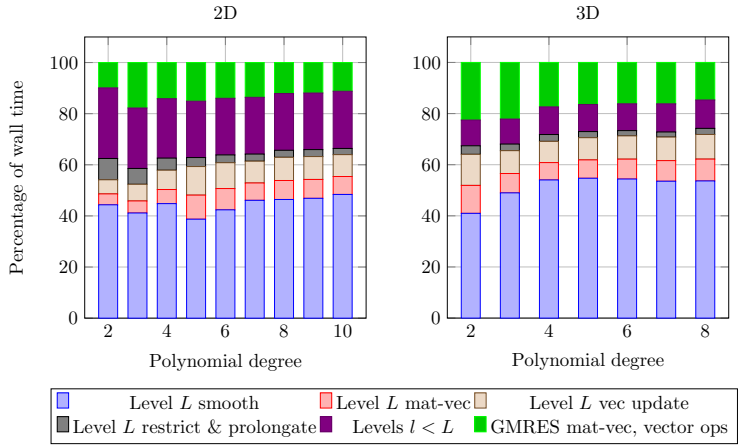}
\caption{Percentage of time spent in various components of the preconditioned GMRES solver with a vertex-patch smoother in 2D and 3D.}
\label{fig:breakdown_GMRES_components_ratio}
\end{figure}

This test uses the \textit{Fused} kernel proposed in the previous section for local smoothing and also serves as a reference implementation for better comparison with the optimisations to be presented afterwards. The operations in the multigrid preconditioner are performed with single precision, while the outer GMRES iterations are in double precision.

First, as depicted in~\cref{fig:breakdown_GMRES_components}, the solver exhibits superior performance when applied to higher-order polynomials in both two and three dimensions. This is because even though higher-order polynomials involve more operations, they need fewer iterations to reach convergence~(\cref{tab:LMG_v_2d,tab:LMG_v_3d}).

Moreover, the smoothing operations contribute substantially to the overall computing time at the finest level, as illustrated in~\cref{fig:breakdown_GMRES_components_ratio}, since the multiplicative vertex patch smoother involves streaming all degrees of freedom four times more than the matrix-vector (mat-vec) operations in two dimensions and eight times more in three dimensions. This stark contrast in computational demands clarifies why smoothing operations occupy a significant portion of the runtime and serve as a primary target for further optimization efforts.

\begin{figure}[tp]
\centering
\includegraphics[width=.45\textwidth]{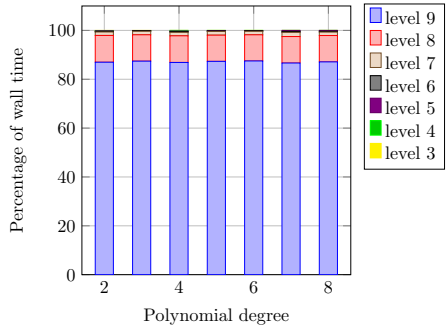}
\caption{Percentage of time spent in various level of the preconditioned GMRES solver with a vertex-patch smoother in 3D.}
\label{fig:breakdown_GMRES_level}
\end{figure}

Additionally, it is evident from~\cref{fig:breakdown_GMRES_level} that operations conducted on the finest level contributes to over 50\% of the total time and 95\% on the last two levels. Operations on coarser levels, on the other hand, are constrained by latency issues since the available work is insufficient to effectively utilize the GPU's capabilities, making these operations better suited for latency-optimized processors such as CPUs.
Incorporating a hybrid implementation that utilizes both CPU and GPU processing could potentially optimize performance~\cite{Sakharnykh16}. However, in optimizing the less than 5\% of the total runtime spent on the coarse girds there is very little to be gained in our case. Therefore, for the present study, a hybrid implementation is not employed.

\subsection{Optimization of smoothing kernels}
\subsubsection{Linear thread indexing}

Divergence caused by conditional statements can have a considerable effect on instruction throughput because it causes threads within the same warp to follow different execution paths. When this divergence occurs, the different paths must be processed consecutively, resulting in an increase in the total number of instructions executed in that warp. 

During the computation of the local residual, all threads are active, while only the interior degrees of freedom are updated during the local solver stage, leading to some threads becoming idle. A simple approach is masking out the threads operating on boundary data, as shown in~\cref{fig:bank_conflict_Q3}a, where the blank area represents the masked threads.

However, this approach presents a critical issue. Within each warp, there are two or more active threads with different branch targets, leading to low branch efficiency. This has been verified by the profiler in~\cref{tab:perf_Q3}. To address this, instead of masking boundary threads, a technique using linear thread indices for indexing interior degrees of freedom is employed. This approach results in the last few warps being idle, with only one warp experiencing divergent branches. As a consequence, the overhead for this method is minimal. This optimization is labeled \emph{Linear Indexing}.
The performance improvement achieved by this approach is noteworthy. The branch efficiency increases from 38\% to 60\%, resulting in a 10\% reduction in runtime, as indicated in~\cref{tab:perf_Q3} with the left two columns.
\begin{table}[tp]
    \centering
    \footnotesize
    \caption{Details of performance analysis for various {\it Fused} kernels with $\mathbb{Q}_3$ element in 3D.}
    \label{tab:perf_Q3}
    \begin{tabular}{lcccc}
    \hline
      & Basic & Linear Indexing & Conflict Free \\ 
    \hline 
    L1 Cache Throughput [\%] & 97.03 & 95.98 & 95.27 \\ 
    FP64 Peak Performance  [\%] & 19 & 21 & 24 \\ 
    Bank Conflicts [per CTA]$^*$ & 113 & 163 & 0 \\ 
    Branch Efficiency [\%] & 38.46  & 63.16 & 60.87 \\ 
    Duration [msecond] & 28.23 & 25.35 & 21.50 \\ 
    \hline
    \multicolumn{4}{l}{\footnotesize * L1 Wavefronts Shared Excessive instruction-level metric.}
    \end{tabular}
\end{table}

\subsubsection{Bank Conflicts}

\begin{figure}[tp]
\centering
    \includegraphics[width=0.9\textwidth]{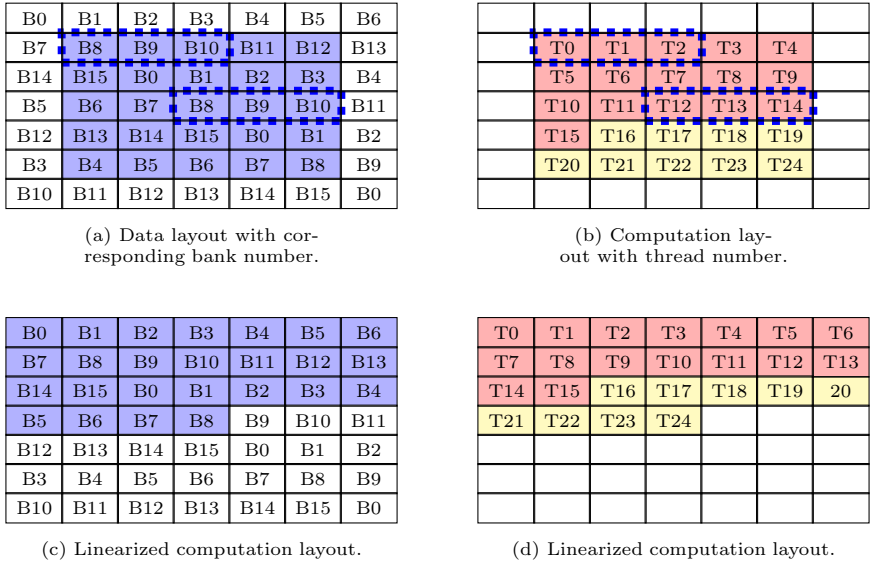}
  \caption{Computation layout and data layout of local solver for $\mathbb{Q}_3$ element.}
  \label{fig:bank_conflict_Q3}
\end{figure}

The access pattern to on-chip memory, also known as \textit{shared memory} in CUDA, plays a crucial role in the overall computation process, given that all calculations occur within this memory. On-chip memory is organized into multiple equally-sized units called \textit{banks}, enabling simultaneous access by threads within a warp. Optimal performance is achieved when all 32 threads access different banks or when they all access the same memory address (broadcast), as these scenarios allow for full-speed memory requests. However, if multiple threads within a warp access the same bank, the accesses are serialized, leading to a phenomenon known as \textit{bank conflict}. The worst case is all threads accessing different addresses on the same bank.

\Cref{fig:bank_conflict_Q3}a illustrates the data layout in on-chip memory and its corresponding bank numbers for the local solver of the $\mathbb{Q}_3$ element. To simplify this discussion, we assume a double-precision data type, which reduces the number of banks from 32 to 16. It is important to note that this optimization also applies to single precision. The values corresponding to the blue part are particularly relevant, as the local solver exclusively affects interior nodes. As discussed in the previous section, linear thread indexing is employed in the computational layout, as shown in~\cref{fig:bank_conflict_Q3}b. Although we illustrate the concept in 2D, it is also applicable to 3D, as the computation only updates one layer at a time.

The current data layout may lead to conflicting on-chip memory accesses. Specifically, when accessing the data for the 25 interior degrees of freedom, an instruction is issued, generating two wavefronts (red and yellow parts in~\cref{fig:bank_conflict_Q3}b) due to the assumed 16 distinct memory banks in on-chip memory. Bank conflicts can only occur within the same wavefront, meaning either the red or yellow wavefront can experience conflicts. As highlighted in the blue dashed box, for example, thread $T0$ and thread $T12$ access the same bank $B8$ but with different memory addresses, leading to a two-way bank conflict. This hardware behavior splits memory requests with conflicts into two independent requests without conflicts, resulting in a reduction in effective bandwidth by a factor of two.

\begin{figure}[tp]
\centering
    \includegraphics[width=0.95\textwidth]{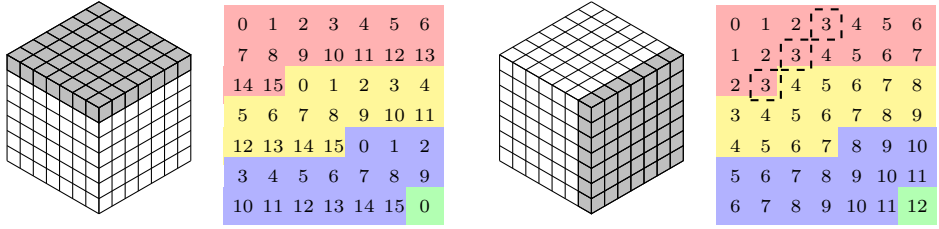}
  \caption{On-chip memory access pattern for tensor contraction with $\mathbb{Q}_3$ element.}
  \label{fig:tensor_mul}
\end{figure}

To achieve a conflict-free access pattern, one can ensure that the interior data remains contiguous in on-chip memory, as demonstrated in~\cref{fig:bank_conflict_Q3}c and~\cref{fig:bank_conflict_Q3}d. This contiguous layout also facilitates coalesced memory access, improving memory access efficiency.
The odd dimension of the matrix facilitates matrix multiplication without encountering bank conflicts, as it ensures that data access does not lead to bank conflicts during the computation. However, bank conflicts can still arise when saving the resulting data back into on-chip memory.

\Cref{fig:tensor_mul} illustrates the on-chip memory access pattern for tensor contractions with the $\mathbb{Q}_3$ element. To maintain brevity, the thread index is omitted, and the numbers indicate the bank numbers. Different colors correspond to different half warps (similar to~\cref{fig:bank_conflict_Q3}), but the two sub-figures are combined into one in this representation. The figure visually demonstrates how the tensor contraction operation may generate bank conflicts during the process of writing the results into on-chip memory.

The desired pattern for tensor contraction is represented on the left in~\cref{fig:tensor_mul}, which avoids bank conflicts. On the other hand, other tensor contraction operations along different dimensions may potentially generate bank conflicts. For instance, the red part in the right plot exhibits a three-way bank conflict.
To address this issue, a general approach is adopted where all tensor contractions are performed using the conflict-free pattern seen on the left. Specifically, during matrix multiplication, the reading of the mass (or stiffness) matrix (\texttt{shape\_data} in Listing \ref{listing:sum_factorized}) is adjusted to ensure conflict-free data access. Since the dimension of the matrix is $(2k+1)$, thus odd and the number of banks is a power of two, reading the mass matrix does not generate bank conflicts. Additionally, when multiple threads request the same address in the same bank, the data is broadcast to the requesting threads.

\begin{figure}[tp]
\centering
\includegraphics[width=.75\textwidth]{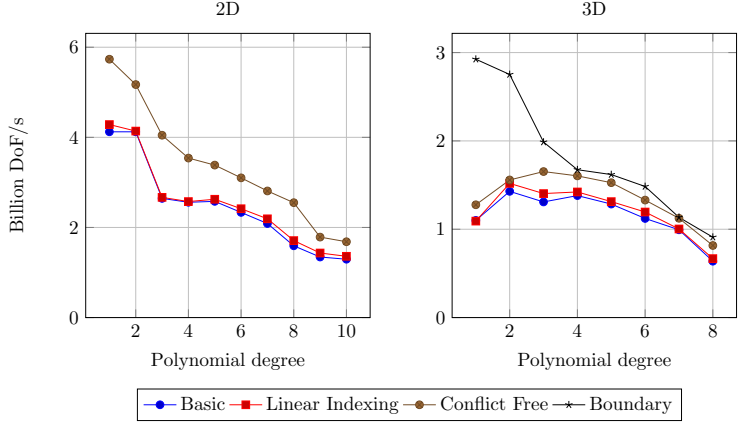}
\caption{Comparison of throughput for different smoother kernels with 26--85 million DoFs in 2D and with 135--721 million DoFs in 3D.}
\label{fig:kernels_throughput}
\end{figure}

To achieve this, sophisticated index calculations are implemented, as shown in Listing \ref{listing:sum_factorized_cf}, which enable conflict-free data access for tensor contraction. This optimization is labeled \emph{Conflict Free}. Moreover, conditional expressions are evaluated at compile-time, ensuring that all threads execute the same instruction without branching.
While matrix operations can be interleaved with vector updates, it is observed that separating the two yields better results. According to~\cref{tab:perf_Q3}, the {\it Conflict Free} kernel is conflict-free, resulting in approximately 15\% reduction in runtime and achieving 24\% of the FP64 peak performance.~\cref{fig:kernels_throughput} provides an overview of the performance results for all polynomial degrees.

\subsubsection{Reduced computation in local residual}
In the case of using an exact local solver, the opportunity arises to further optimize the algorithm and improve its efficiency.
Consider the matrix form of a local operator $\overline A_j$:
\begin{equation*}
    \overline A_j = \begin{pmatrix}
    A^{II} & A^{IB} \\ A^{BI} & A^{BB}
    \end{pmatrix},
\end{equation*}
where $I$ and $B$ refer to the interior and boundary, respectively. Then, the local smoothing step can be rewritten as follows:
\begin{align*}
 x_j + R_j^T A_j^{-1} R_j (b_j - \overline A_j x_j)
= \begin{pmatrix}
    A_j^{-1}(b^I- A^{IB}x^B) \\ x^B
   \end{pmatrix}
\end{align*}
This optimization enables reduced DRAM data loading by eliminating the need for $x^I$ values, resulting in significant operational savings during residual computation. The local operator $A^{IB}$ can still be represented as the tensor product of 1D matrices
\begin{equation}\label{eq:A_ib_2d}
    A^{IB} = (A_1^{II}\otimes M_0^{IB}, A_1^{IB} \otimes M_0^I) + 
                        (M_1^{II}\otimes A_0^{IB}, M_1^{IB} \otimes A_0^I)
\end{equation}
in two dimensions and 
\begin{equation}\label{eq:A_ib_3d}
\begin{aligned}
    A^{IB} = & (A_2^{II}\otimes M_1^{II}\otimes M_0^{IB}, A_2^{II}\otimes M_1^{IB}\otimes M_0^{I}, A_2^{IB}\otimes M_1^{I}\otimes M_0^{I}) + \\
    & (M_2^{II}\otimes A_1^{II}\otimes M_0^{IB}, M_2^{II}\otimes A_1^{IB}\otimes M_0^{I}, M_2^{IB}\otimes A_1^{I}\otimes M_0^{I}) + \\
    & (M_2^{II}\otimes M_1^{II}\otimes A_0^{IB}, M_2^{II}\otimes M_1^{IB}\otimes A_0^{I}, M_2^{IB}\otimes M_1^{I}\otimes A_0^{I})
\end{aligned}
\end{equation}
in three dimensions, where $A^{I} = (A^{II}, A^{IB})$ and $M^{I} = (M^{II}, M^{IB})$. This optimization is labeled \emph{Boundary}.

\begin{figure}[tp]
\centering
\includegraphics[width=.65\textwidth]{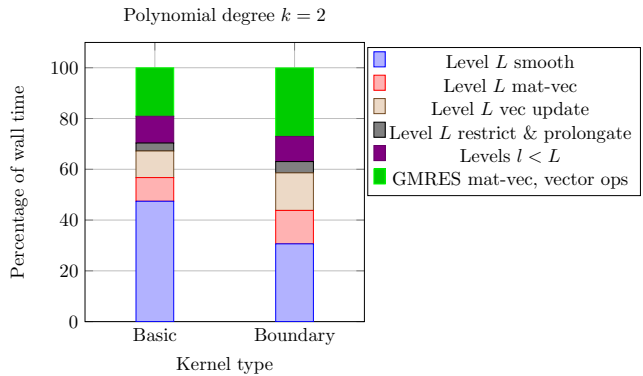}
\caption{Comparison of percentage of time spent in various components of the preconditioned GMRES solver using \textit{Basic} and \textit{Boundary} kernel in 3D.}
\label{fig:breakdown_GMRES_components_ratio_q2}
\end{figure}

However, matrix-free evaluation of operator $A_{IB}$ is more complex than $A_\ell$ for the entire patch. Inconsistent matrix dimensions and increased synchronization make this approach more suitable for lower order elements and less effective for higher-orders, as depicted in~\cref{fig:kernels_throughput} ({\it Boundary}). The feasibility of this optimization also relies on the assumption that $A_j^{-1}$ is accurate, which might not hold true for more complicated problems.
In~\cref{fig:breakdown_GMRES_components_ratio_q2}, we we present a comparison of percentage of different components following optimization. The optimized smoothing operation, once the sole decisive factor in the overall solution time, is no longer dominant. Given that the smoother is not the only bottleneck, it is essential to also explore the optimization of the finite element operator, a primary focus of our upcoming research. 

\begin{figure}[tp]
\centering
\includegraphics[width=.45\textwidth]{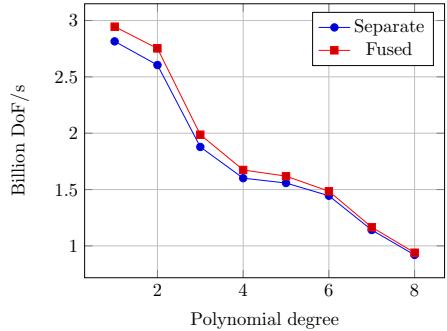}
\caption{Comparison of throughput for optimized {\it Separate} and {\it Fused} kernels with 135–721 million DoFs in 3D.}
\label{fig:sep_fused}
\end{figure}

Finally, \cref{fig:sep_fused} gives comparative results for the \emph{Separated} kernel using the above optimizations, where it can be seen that the optimizations have a more pronounced impact on the reuse of local data, which reduces the accesses to the global vectors, making the \emph{Fused} kernel outperform the separated kernel to all polynomial degrees. 
However, the expected significant performance gain resulting from the reduced global vector accesses is not evident in this scenario. As depicted in~\cref{tab:perf_Q2}, when kernels only involve data access (referred to as Data only), the \emph{Fused} kernel expends, as expected, approximately two-thirds of the time compared to the \emph{Separate} kernel. This time difference theoretically represents the portion of time expected to be saved when computation is integrated, constituting roughly 15\% of the total time. However, given our compute-bound nature, this improvement is not immediately apparent. Additionally, in the case of the local solver, the \emph{Separate} kernel demonstrates the ability to allocate more blocks to Streaming Multiprocessors (SMs) due to its utilization of fewer on-chip memories, further diminishing the disparity between the two approaches.

\begin{table}[tp]
    \centering
    \footnotesize
    \caption{Details of performance analysis for {\it Fused} and {\it Separate} kernel with $\mathbb{Q}_2$ and $\mathbb{Q}_6$ element in 3D.}
    \label{tab:perf_Q2}
    \begin{tabular}{clccc}
    \hline
      & total time [ms] & data only time [ms] & block/SM$^*$ \\
    \hline 
   & Fused & 6.85 & 1.89 & 3 \\
  $\mathbb{Q}_2$ & Separate (residual) & 5.58 & 2.02 & 3 \\
   & Separate (solver) & 1.66 & 1.07 & 8 \\
    \hline
   & Fused & 43.20 & 7.56 & 3\\
  $\mathbb{Q}_6$ & Separate (residual) & 30.20 & 8.28 &  3 \\
   & Separate (solver) & 14.24 & 6.78 & 4 \\
    \hline
    \multicolumn{4}{l}{\footnotesize * Occupancy limit due to on-chip memory usage.}
    \end{tabular}
\end{table}

\subsection{Exploiting mixed precision}
When solving linear systems of equations, whether through iterative methods or direct methods, maintaining high floating-point precision is crucial to obtain accurate results, as rounding errors can accumulate during computations. G{\"o}ddeke~\cite{goddeke2007performance} compares native double-precision algorithms with mixed-precision algorithms for solving linear systems using CG and multigrid methods with a Jacobi smoother. The proposed mixed precision approach works very well while maintaining the same accuracy as the reference solver running in double precision.

\begin{table}[tp]
    \centering
    \footnotesize
    \caption{Comparison of Mixed and Double Precision for full multigrid method in three dimensions.}
    \label{tab:mixed_double}
    \begin{tabular}{lcccccccc}
    \hline
          &  & \multicolumn{3}{c}{ Double Precision } & \multicolumn{3}{c}{ Mixed Precision } &  \\
          degree & DoF & $L_2$ error & iteration & time[s]  & $L_2$ error & iteration & time[s] & Speedup\\
    \hline
         1 & 135M & 1.12e-06 & 5 & 3.391 & 1.12e-06 & 5 & 2.385 & 1.42 \\
         3 & 454M & 2.30e-13 & 3 & 3.941 & 2.30e-13 & 3 & 2.418 & 1.59 \\
         7 & 721M & 2.89e-16 & 2 & 5.891 & 2.89e-16 & 2 & 3.326 & 1.77 \\   
    \hline
    \end{tabular}
\end{table}
In~\cref{tab:mixed_double}, we compare a GMRES method run in double precision preconditioned by a multigrid V-cycle with one pre- and postsmoothing step in single precision against doing all computations in the multigrid V-cycle in double precision. The right hand side of the test problem is manufactured such that the exact solution $u$ is $u(\mathbf{x}) = \prod_{i = 1}^{d} \mathrm{sin}(\pi x_i)$. A relative residual reduction of $10^{-9}$ uses as the stopping criterion.
We find that the final $L_2$ discretization error remains the same, regardless of whether the V-cycle is done with single or double precision. This suggests that the mixed-precision approach is capable of achieving the same level of accuracy as the double-precision approach. Moreover, the mixed-precision approach offers significant advantages in terms of runtime. Specifically, it can speed up the solving of linear systems by 42\% to 77\% compared to using double precision throughout the entire multigrid V-cycle.

\begin{figure}[tp]
\centering
\includegraphics[width=.75\textwidth]{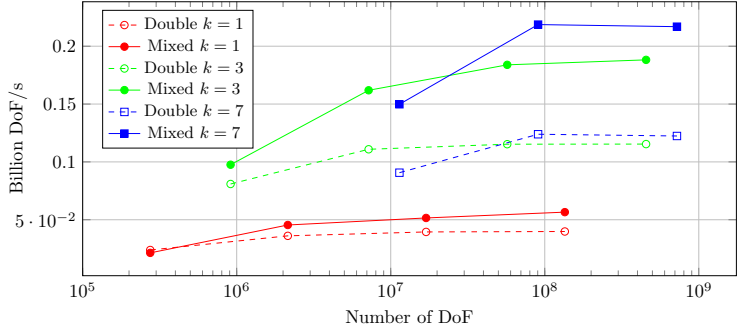}
\caption{Comparison of Double and Mixed Precision performance for solving the Poisson equation with preconditioned GMRES method in three dimensions.}
\label{fig:throughput}
\end{figure}

\Cref{fig:throughput} presents the performance results for solving the Poisson equation, measured in seconds per billions of degrees of freedom. For all other polynomial degrees, the performance is similar to that of $k=7$.
Both double and mixed precision solutions exhibit relatively longer times for small problem sizes, where there are not enough parallelism in the problem to fully saturate the hardware. For higher-order elements, employing mixed precision can reduce the solution time by nearly a factor of two. The performance curve increases smoothly, benefiting from the well-implemented solution and the flexibility of CUDA, which ensures that the performance is not affected by cache size limitations as seen in CPUs. This consistent performance holds across all order elements. Although local smoothing operations are more computationally expensive for higher-order elements, the benefit lies in the fact that fewer iterations are required to achieve convergence, compensating for the increased computational cost.

\subsection{Roofline analysis}
\begin{figure}[tp]
\centering
\includegraphics[width=.45\textwidth]{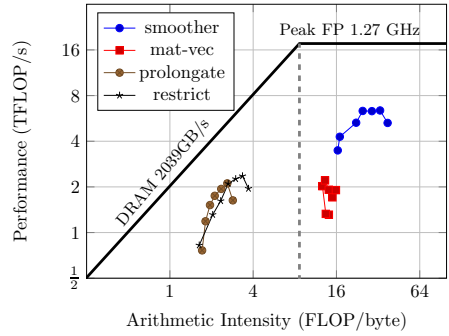}
\caption{Roofline model for different operations of the 3D Laplacian on Ampere A100 GPU at 1.27 GHz, displaying the data for polynomial degrees $k = 2,3,...,8$. The problem sizes are between 135 and 721 million DoFs.}
\label{fig:roofline_3d}
\end{figure}

\Cref{fig:roofline_3d} displays the performance of operations in multigrid method in terms of the roofline performance model. By exploiting the tensor product structure, the matrix-free multiplication operation and smoothing operation are not limited by DRAM bandwidth. Again, because of the heavy use of on-chip memory as indicated in~\cref{tab:perf_Q3}, we believe that on-chip memory bandwidth is the real limiting factor. The conventional roofline model is under the assumption that the DRAM bandwidth is the limiting factor. 

\begin{figure}[tp]
\centering
\includegraphics[width=.45\textwidth]{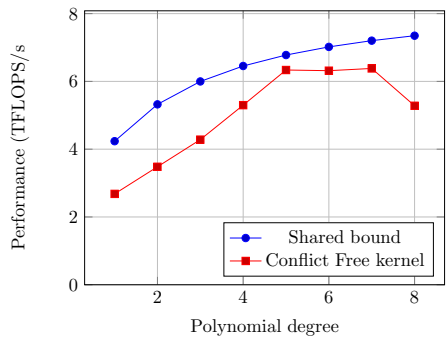}
\caption{On-chip memory roofline model and achieved performance for smoothing kernel in 3D, displaying the data for polynomial degrees $k = 1, 2,...,8$. The problem sizes are between 135 and 721 million DoFs.}
\label{fig:roofline_shmem}
\end{figure}

Inspired by~\cite{swirydowicz2019acceleration}, we introduce a roofline model that specifically focuses on the on-chip memory bandwidth, enabling a more accurate understanding of the factors affecting the performance.
The bandwidth of on-chip memory can be estimated as follows:
\begin{equation*}
    B = \# \text{SMs} \times \# \text{banks} \times \text{word length} \times \text{clock speed}.
\end{equation*}
For Ampere A100 GPU, the corresponding bandwidth is $108\cdot 32 \cdot 4 \cdot 1.27=17.145 \text{TB/s}$, which aligns with the measured on-chip bandwidth of 127.7 bytes/clk/SM as reported in~\cite{sun2022dissecting}. Then, the on-chip memory roofline model is given by:
\begin{equation}
    \mathcal{R} = B \cdot \frac{F}{d_r+d_w},
\end{equation}
Where $d_r$ denotes the number of bytes read from on-chip memory, $d_w$ denotes the number of bytes written to on-chip memory, and $F$ denotes the number of floating point operations which generate the on-chip traffic of $d_r+d_w$ bytes. \Cref{fig:roofline_shmem} shows the roofline of the \textit{Conflict Free} kernel, from which we see that the performance results are very close to the roofline model for higher order polynomial degrees $k=5, 6, 7$. Also the profiler shows that the L1/TEX cache throughput is more than 95\%, thereby affirming on-chip memory bandwidth as the predominant constraining factor. Optimal performance reaches a notable 36\% of peak performance holding true for both double and single precision. A decline in performance occurs when polynomial degree $k=8$. This is attributed to the excessive utilization of on-chip memory, leading to a substantial reduction in the concurrent processing capacity of blocks per multiprocessor.

\section{Conclusion}
\label{sec:conclusion}
In this work we have presented a GPU implementation of tensor product based vertex patch smoothers for higher order finite element methods in two and three dimensions on equidistant Cartesian meshes. We have discussed various optimization strategies to improve the performance of the methods, ranging from numerical techniques  to programming considerations within the CUDA framework. First, efficient local solver through fast diagonalization minimizes computational and space complexity, resulting in a significant gain in algorithm efficiency, which makes it feasible for high order elements in three dimensions. Then, a colorized algorithm with localized operations, yielding remarkable speedups of up to 2.6x in 2D and 7x in 3D. 

Detailed performance analysis indicates that the obtained method preforms equally well for high order polynomial degrees. Furthermore, an optimal on-chip memory access pattern with zero bank conflict enables us utilize up to 36\% of the theoretical FP32/FP64 peak performance. Using the mixed precision approach, a speedup of up to 77\% is achieved in solving the 3D Poisson problem with 721M DoF, taking only 3.3 seconds. In addition, by exploiting the specificity of the exact local solver, fewer operations and data are required to compute the local residuals, making the method also attractive for lower order elements.

\bibliographystyle{siam}
\bibliography{references}

\appendix
\section{Device code for tensor contraction of local smoothing operation in 3D}
In Listing~\ref{listing:sum_factorized_cf}, we offer partial code examples that demonstrate how to eliminate bank conflicts by reordering local operations for improved performance with the straightforward implementation Listing~\ref{listing:sum_factorized} as a comparison.
\footnotesize
\begin{lstlisting}[caption={Device code for tensor contraction operations in 3D (Basic).},label={listing:sum_factorized}]
template <typename Number, int n_dofs_1d, int dir, bool add, bool sub>
__device__ void apply(const Number *shape_data, const Number *in, Number *out) {
    ...
    Number pval[n_dofs_1d] = {};
    for (unsigned int z = 0; z < n_dofs_1d; ++z) 
        for (unsigned int k = 0; k < n_dofs_1d; ++k) {
            const unsigned int shape_idx = row * n_dofs_1d + k;
            const unsigned int source_idx =
                (dir == 0) ? (col * n_dofs_1d + k + z * stride) :
                (dir == 1) ? (k * n_dofs_1d + col + z * stride) :
                             (z * n_dofs_1d + col + k * stride);
            pval[z] += shape_data[shape_idx] * in[source_idx];
        }
    for (unsigned int z = 0; z < n_dofs_1d; ++z) {
        const unsigned int destination_idx =
            (dir == 0) ? (col * n_dofs_1d + row + z * stride) :
            (dir == 1) ? (row * n_dofs_1d + col + z * stride) :
                         (z * n_dofs_1d + col + row * stride);
        if (add) out[destination_idx] += pval[z];
        else if (sub) out[destination_idx] -= pval[z];
        else out[destination_idx] = pval[z];
    }
}
\end{lstlisting}
\begin{lstlisting}[caption={Device code for tensor contraction operations in 3D (Conflict-Free).},label={listing:sum_factorized_cf}]
template <typename Number, int n_dofs_1d, int dir, bool add, bool sub>
__device__ void apply(const Number *shape_data, const Number *in, Number *out) {
    ...
    Number pval[n_dofs_1d] = {};
    for (unsigned int z = 0; z < n_dofs_1d; ++z)
        for (unsigned int k = 0; k < n_dofs_1d; ++k) {
            const unsigned int shape_idx =
                (dir == 0) ? (col * n_dofs_1d + k) :
                (dir == 1) ? (row * n_dofs_1d + k) :
                             (z * n_dofs_1d + k);
            const unsigned int source_idx =
                (dir == 0) ? (row * n_dofs_1d + k + z * stride) :
                (dir == 1) ? (k * n_dofs_1d + col + z * stride) :
                             (row * n_dofs_1d + col + k * stride);
            pval[z] += shape_data[shape_idx] * in[source_idx];
        }
    for (unsigned int z = 0; z < n_dofs_1d; ++z) {
        const unsigned int destination_idx = 
            row * n_dofs_1d + col + z * stride;
        if (add) out[destination_idx] += pval[z];
        else if (sub) out[destination_idx] -= pval[z];
        else out[destination_idx] = pval[z];
    }
}
\end{lstlisting}

\end{document}